# LIPSCHITZ SPACES AND CALDERÓN-ZYGMUND OPERATORS ASSOCIATED TO NON–DOUBLING MEASURES

JOSÉ GARCÍA-CUERVA AND A. EDUARDO GATTO

ABSTRACT. In the setting of $\mathbb{R}^d$ with an $n-$dimensional measure $\mu$, we give several characterizations of Lipschitz spaces in terms of mean oscillations involving $\mu$. We also show that Lipschitz spaces are preserved by those Calderón-Zygmund operators $T$ associated to the measure $\mu$ for which $T(1)$ is the Lipschitz class 0.

## 1. INTRODUCTION

In the sequel $\mathbb{R}^d$ will denote the euclidean space of vectors $x = (x_1, \ldots, x_d)$ with $d$ real components and $\mu$ a Radon $n-$dimensional measure on $\mathbb{R}^d$. We say that $\mu$ is an $n-$dimensional measure, $n \in \mathbb{R}$, $0 < n \leq d$, if it satisfies the following growth condition: There is a constant $C$ such that

$$(1.1) \qquad \mu(B(x,r)) \leq Cr^n$$

for every ball $B(x,r)$ with center $x \in \mathbb{R}^d$ and radius $r > 0$. This allows, in particular, non-doubling measures.

The study of Calderón-Zygmund operators associated to an $n-$dimensional measure was carried out, in the Lebesgue spaces, by Nazarov, Treil and Volberg (see [NTV1, NTV2]) and also by Tolsa (see [To1, To2]). Further results, dealing with $BMO$ and $H^1$ and providing boundedness criteria in the spirit of the $T(1)$ or $T(b)$ theorems, were obtained as well (see [NTV3, MMNO, To3])

*Date*: November 1, 2018.

1991 *Mathematics Subject Classification.* 2000 revision, 42B20, 26B35, 47B38, 47G10.

Supported in part by DGES, Spain, under grant BFM2001-0189. It is a pleasure for the second author, to thank the members of the Mathematics Department of Universidad Autónoma de Madrid for their friendly hospitality.

In [GG] we have also studied the theory of fractional integral operators associated to an $n-$dimensional measure $\mu$ on Lebesgue spaces and Lipschitz spaces.

The present paper is devoted to further study the Lipschitz spaces associated to an $n-$dimensional measure $\mu$ on $\mathbb{R}^d$ and Calderón-Zygmund operators associated to $\mu$ on these spaces.

In section 2, we obtain several characterizations in terms of mean oscillations of functions in the Lipschitz spaces associated to $\mu$. Then these characterizations are used in section 3 to establish the boundedness of Calderón-Zygmund operators.

From now on, all balls that we consider will be centered at points in the support of $\mu$.

## 2. Characterization of Lipschitz spaces

In order to prove the main theorem of this section we will need the following known definition and lemma (see [To3])

**Definition 2.1.** Let $\beta$ be a fixed constant. A ball $B$ is called $\beta-$doubling if
$$\mu(2B) \leq \beta\mu(B).$$

**Lemma 2.2.** *Let $f \in L^1_{\mathrm{loc}}(\mu)$. If $\beta > 2^d$, then, for almost every $x$ with respect to $\mu$, there exists a sequence of $\beta-$doubling balls $B_j = B(x, r_j)$ with $r_j \to 0$, such that*
$$\lim_{j\to\infty} \frac{1}{\mu(B_j)} \int_{B_j} f(y) \,\mathrm{d}\mu(y) = f(x).$$

**Proof.** We will show that for almost every $x$ with respect to $\mu$ there is a $\beta-$doubling ball centered at $x$ with radius as small as we wish. This fact, combined with the differentiation theorem, completes the proof of the lemma.

We know that for almost every $x$ with respect to $\mu$

(2.1) $$\lim_{r\to 0} \frac{\mu(B(x,r))}{r^d} > 0$$

(differentiation of $\mu$ with respect to Lebesgue measure, see [M]). Now for $x$ satisfying (2.1) take $B = B(x, r)$ and assume that none of the balls $2^{-k}B$, $k \geq 1$, is $\beta-$doubling. Then it easy to see that $\mu(B) > \beta^k \mu(2^{-k}B)$ for all $k \geq 1$. Therefore
$$\frac{\mu(2^{-k}B)}{(2^{-k}r)^d} < \left(\frac{2^d}{\beta}\right)^k \frac{\mu(B)}{r^d}.$$

Note that, since $\beta > 2^d$, the right hand side tends to zero for $k \to \infty$, which is a contradiction.

□

Now we can state and prove the main result of this section.

**Theorem 2.3.** *For a function $f \in L^1_{\text{loc}}(\mu)$, the conditions I, II, and III below, are equivalent*

(**I**) *There exist some constant $C_1$ and a collection of numbers $f_B$, one for each ball $B$, such that these two properties hold: For any ball $B$ with radius $r$*

$$\frac{1}{\mu(2B)} \int_B |f(x) - f_B| \, d\mu(x) \leq C_1 r^\alpha, \tag{2.2}$$

*and for any ball $U$ such that $B \subset U$ and $\text{radius}(U) \leq 2r$.*

$$|f_B - f_U| \leq C_1 r^\alpha, \tag{2.3}$$

(**II**) *There is a constant $C_2$ such that*

$$|f(x) - f(y)| \leq C_2 |x - y|^\alpha \tag{2.4}$$

*for $\mu$-almost every $x$ and $y$ in the support of $\mu$.*

(**III**) *For any given $p$, $1 \leq p \leq \infty$, there is a constant $C(p)$, such that for every ball $B$ of radius $r$, we have*

$$\left(\frac{1}{\mu(B)} \int_B |f(x) - m_B(f)|^p \, d\mu(x)\right)^{1/p} \leq C(p) r^\alpha, \tag{2.5}$$

*where $m_B(f) = \dfrac{1}{\mu(B)} \int_B f(y) \, d\mu(y)$ and also for any ball $U$ such that $B \subset U$ and $\text{radius}(U) \leq 2r$.*

$$|m_B(f) - m_U(f)| \leq C(p) r^\alpha, \tag{2.6}$$

*In addition, the quantities: $\inf C_1$, $\inf C_2$, and $\inf C(p)$ with a fixed $p$ are equivalent.*

**Proof.** (I) ⇒ (II). Consider $x$ as in the Lemma and let $B_j = B(x, r_j)$, $j \geq 1$, a sequence of $\beta$-doubling balls with $r_j \to 0$. We will show first that (2.2) implies

$$lim_{j \to \infty} f_{B_j} = f(x).$$

It suffices to observe that

$$|m_{B_j}(f) - f_{B_j}| \leq \frac{1}{\mu(B_j)} \int_{B_j} |f(y) - f_{B_j}| \, d\mu(y)$$

$$\leq \frac{\mu(2B_j)}{\mu(B_j)} \frac{1}{\mu(2B_j)} \int_{B_j} |f(y) - f_{B_j}| \, d\mu(y) \leq \beta C_1 r_j^\alpha.$$

Next, let $x$ and $y$ be two points as in the lemma. Take $B = B(x, r)$ any ball with $r \leq |x - y|$ and let $U = B(x, 2|x - y|)$. Now define $B_k = B(x, 2^k r)$, for $0 \leq k \leq \bar{k}$, where $\bar{k}$ is the first integer such that $2^{\bar{k}} r \geq |x - y|$. Then

$$|f_B - f_U| \leq \sum_{k=0}^{\bar{k}-1} |f_{B_k} - f_{B_{k+1}}| + |f_{B_{\bar{k}}} - f_U|$$

$$\leq C_1 \sum_{k=0}^{\bar{k}} (2^k r)^\alpha \leq C' C_1 |x - y|^\alpha,$$

where $C'$ is independent of $x$ and $B$.

A similar argument can be made for the point $y$ with any ball $B' = B(y, s)$ such that $s \leq |x - y|$ and $V = B(y, 3|x - y|)$. Therefore

$$|f_B - f_{B'}| \leq |f_B - f_U| + |f_U - f_V| + |f_V - f_{B'}| \leq C'' C_1 |x - y|^\alpha.$$

Finally, take two sequences of $\beta$-doubling balls $B_j = B(x, r_j)$ and $B'_j = B(y, s_j)$ with $r_j \to 0$ and $s_j \to 0$. We have

$$|f(x) - f(y)| = \lim_{j \to \infty} \left| f_{B_j} - f_{B'_j} \right| \leq C'' C_1 |x - y|^\alpha.$$

(II) $\Rightarrow$ (III). It is immediate. Note also that (II) $\Rightarrow$ (I) is immediate as well.

(III) $\Rightarrow$ (I). Define first $f_B = m_B(f)$. Then (2.3) is exactly (2.6). To prove (2.2), we write

$$\frac{1}{\mu(2B)} \int_B |f(y) - f_B| \, d\mu(y) \leq \frac{1}{\mu(2B)} \left( \int_B |f(y) - f_B|^p \right)^{1/p} \mu(B)^{1/p'}$$

$$\leq \frac{\mu(B)}{\mu(2B)} \left( \frac{1}{\mu(B)} \int_B |f(y) - f_B|^p \, d\mu(y) \right)^{1/p} \leq C(p) r^\alpha.$$

This concludes the proof of the theorem.

$\square$

**Remark 2.4.** Theorem 2.3 is also true if the number 2 in condition (I) is replaced by any fixed $\rho > 1$. In that case, the proof uses $(\rho, \beta)$-doubling balls, that is, balls satisfying $\mu(\rho B) \leq \beta \mu(B)$. However this extension is not needed in our paper.

The idea of combining the mean oscillation condition with an extra condition as in (I) originates in the work of Tolsa [To3] on regular $BMO$, whereas the introduction of the $\rho$ factor in (I) comes from [NTV3].

**Definition 2.5.** We shall call Lipschitz function of order $\alpha$ with respect to $\mu$ to a function, or rather the corresponding Lebesgue class in $L^1_{\text{loc}}(\mu)$, which satisfies any, and hence all, of the conditions of theorem 2.3.

The linear space of all Lipschitz functions of order $\alpha$, with respect to $\mu$, modulo constants, becomes, with the norm $\inf C_2$ of Theorem 2.3, a Banach space, which we shall call $\mathcal{L}ip(\alpha, \mu)$.

**Remark 2.6.** It is easy to see that $\mathcal{L}ip(\alpha, \mu)$ coincides with the space of Lipschitz functions of order $\alpha$ on the support of $\mu$.

Note that by the extension theorem of Banach (see [B] or [Mi]), any Lipschitz function of order $\alpha$ with respect to $\mu$ has an extension to $\mathbb{R}^d$ that is a Lipschitz function of order $\alpha$ with an equivalent norm.

**Remark 2.7.** For $0 < \alpha \leq 1$, a telescoping argument like the one used in the proof of (I) $\Rightarrow$ (II) in Theorem 2.3 shows that (2.3) is equivalent to

(2.7) $$|f_B - f_U| \leq C'_1 \text{radius}(U)^\alpha$$

for any two balls $B \subset U$. Furthermore (2.7) is also equivalent to

(2.8) $$|f_B - f_U| \leq C''_1 K_{B,U} \text{radius}(U)^\alpha,$$

for any two balls $B \subset U$, where $K_{B,U}$ is the constant introduced by X. Tolsa in [To3], given by

$$K_{B,U} = 1 + \sum_{j=1}^{N_{B,U}} \frac{\mu(2^j B)}{(2^j B)^n},$$

with $N_{B,U}$ equal to the first integer $k$ such $2^k \text{radius}(B) \geq \text{radius}(U)$. Indeed (2.8) for comparable balls, that is, for $\text{radius}(U) \leq 2\text{radius}(B)$, reduces to (2.3) because, in this case, $K_{B,U}$ is controlled by an absolute constant.

Note that (2.2) and (2.8) make sense also for $\alpha = 0$ and the space defined by them is the space $RBMO(\mu)$ of X. Tolsa (see [To3]). Therefore, the spaces $\mathcal{L}ip(\alpha, \mu)$, $0 < \alpha \leq 1$ can be seen as members of a family containing also $RBMO(\mu)$.

## 3. Calderón-Zygmund operators

**Definition 3.1.** An "$n$-dimensional" standard kernel $K$ on $\mathbb{R}^d$ will be a function $K \in L^1_{\text{loc}}(\mathbb{R}^d \times \mathbb{R}^d \setminus \{(x, y) : x = y\}, \mu \times \mu)$, satisfying the following two conditions for some constants $C > 0$ and $\varepsilon \in ]0, 1]$:

$$(i)\ |K(x,y)| \leq \frac{C}{|x-y|^n}, \quad \text{for all } x \neq y \text{ and}$$

$$(ii)\ |K(x,y) - K(x',y)| \leq C\frac{|x-x'|^\varepsilon}{|x-y|^{n+\varepsilon}} \quad \text{for } |x-y| \geq 2|x-x'|.$$

When we want to be precise about $\varepsilon$, we shall say that $K$ is an $n-$dimensional, $\varepsilon-$regular, standard kernel.

**Definition 3.2.** By a Calderón-Zygmund operator on $\mathbb{R}^d$ with respect to $\mu$, we shall mean a linear operator $T$, bounded on $L^2(\mu)$, which is associated to an $n-$dimensional standard kernel $K$ in the sense that for every $f \in L^2(\mu)$, with compact support

$$Tf(x) = \int_{\mathbb{R}^d} K(x,y)\, f(y)\, d\mu(y), \quad \text{for } \mu\text{-almost every } x \in \mathbb{R}^d \setminus \operatorname{supp} f.$$

**Remark 3.3.** We have not required $\epsilon-$regularity of $K$ in the second variable, as is usually done, because we do not need it for proving Theorem 3.6.

The definition of a Calderón-Zygmund operator $T$ can be extended to Lipschitz functions in the following way:

**Definition 3.4.** Assume that $T$ is a Calderón-Zygmund operator associated to $\mu$ having an $\epsilon-$regular kernel. First, for each ball $B$, with center $x_B$, $x \in B$ and $f \in \mathcal{L}ip(\alpha,\mu)$, $0 < \alpha < \epsilon$, we define

$$(3.1)\quad T_B(f)(x)$$
$$= T(f\chi_{2B})(x) + \int_{\mathbb{R}^d \setminus 2B} \left(K(x,z) - K(x_B,z)\right) f(z)\, \mathrm{d}\mu(z)$$

Now $T(f)$ is defined as the class of all functions $g$ on $\mathbb{R}^d$ such that for any ball $B$ there is a constant $c_B$ satisfying

$$g(x) = T_B(f)(x) + c_B \quad \mu - \text{almost everywhere on } B.$$

**Remark 3.5.** In order to see that $T(f)$ is well defined, it suffices to realize that if $B \subset V$ are balls, then $T_B(f) - T_V(f)$ is constant $\mu-$a. e. over $B$. Then, considering an increasing sequence of balls whose union is $\mathbb{R}^d$, we can see that there is a function $g$ such that $g - T_B(f)$ is constant $\mu-$a. e. on $B$ for each ball $B$ in the sequence. Consider now any ball $B$. Since this ball $B$ is contained in some ball of the sequence, by our previous observation, $g - T_B(f)$ is constant $\mu-$a. e. on $B$.

**Theorem 3.6.** *Let $T$ be a Calderón-Zygmund operator associated to $\mu$ having an $\epsilon-$regular kernel and $0 < \alpha < \epsilon \leq 1$. Then $T$ is bounded on $\mathcal{L}ip(\alpha,\mu)$ if and only if $T(1) = 0$ in $\mathcal{L}ip(\alpha,\mu)$.*

**Proof.** The fact that $T(1) = 0$ in $\mathcal{L}ip(\alpha, \mu)$ is a necessary condition follows immediately from the fact that $\|1\|_{\mathcal{L}ip(\alpha,\mu)} = 0$.

To prove that this condition is sufficient, take $f \in \mathcal{L}ip(\alpha, \mu)$ and two balls $B \subset U$ with radius$(U) \leq 2r$, where $r$ is the radius of $B$. We want to show that there are two constants $a_B$ and $a_U$ such that

$$(3.2) \qquad \frac{1}{\mu(2B)} \int_B |T_U(f)(x) - a_B| \, d\mu(x) \leq C \, \|f\|_{\mathcal{L}ip(\alpha,\mu)} \, r^\alpha$$

and

$$(3.3) \qquad |a_B - a_U| \leq C \, \|f\|_{\mathcal{L}ip(\alpha,\mu)} \, r^\alpha$$

with $C$ independent of $B, U$ and $f$.

Since $T(f)(x) = T_U(f)(x) + c_U$ on $U$, where $c_U$ is a constant, then (3.2) and (3.3) follow for $T(f)$ with constants $a'_B = a_B + c_U$ and $a'_U = a_U + c_U$.

By the hypothesis $T(1) = 0$, we can assume that $\int_{2U} f \, d\mu = 0$. Let's now define

$$a_B = \frac{1}{\mu(2B)} \int_B T_U(f)(x) \, d\mu(x)$$

and, similarly $a_U$ with $U$ instead of $B$. Observe that

$$\frac{1}{\mu(2B)} \int_B |T_U(f)(x) - a_B| \, d\mu(x)$$
$$\leq \frac{1}{\mu(2B)} \int_B |T_U(f)(x)| \, d\mu(x) + \frac{1}{\mu(2B)} \int_B |a_B| \, d\mu(x)$$
$$\leq 2 \frac{1}{\mu(2B)} \int_B |T_U(f)(x)| \, d\mu(x),$$

since $\frac{\mu(B)}{\mu(2B)} \leq 1$. On the other hand

$$|a_B - a_U| \leq \frac{1}{\mu(2B)} \int_B |T_U(f)(x)| \, d\mu(x) + \frac{1}{\mu(2U)} \int_U |T_U(f)(x)| \, d\mu(x).$$

Now we want to estimate the quantity

$$(3.4) \qquad \frac{1}{\mu(2B)} \int_B |T_U(f)(x)| \, d\mu(x),$$

where $B$ can also be $U$. To do that we write

$$T_U(f)(x) = T(f\chi_{2B})(x) + T(f\chi_{2U \setminus 2B})(x)$$
$$+ \int_{\mathbb{R}^d \setminus 2U} (k(x, z) - k(x_U, z)) \, f(z) \, d\mu(z).$$

We then have

$$\frac{1}{\mu(2B)} \int_B |T_U(f)(x)| \, d\mu(x)$$
$$\leq \frac{1}{\mu(2B)} \int_B |T(f\chi_{2B})(x)| \, d\mu(x) + \frac{1}{\mu(2B)} \int_B |T(f\chi_{2U\setminus 2B})(x)| \, d\mu(x)$$
$$+ \frac{1}{\mu(2B)} \int_B \left| \int_{\mathbb{R}^d\setminus 2U} (k(x,z) - k(x_U, z)) f(z) \, d\mu(z) \right| d\mu(x) = A_1 + A_2 + A_3.$$

To estimate $A_1$, we use that $T$ is bounded in $L^2(\mu)$ and the fact that $|m_{2B}(f)| \leq C \|f\|_{\mathcal{L}ip(\alpha,\mu)} r^\alpha$, which follows from $|m_{2B}(f) - m_{2U}(f)| \leq C \|f\|_{\mathcal{L}ip(\alpha,\mu)} r^\alpha$ and $m_{2U}(f) = 0$. We have

$$\frac{1}{\mu(2B)} \int_B |T(f\chi_{2B})(x)| \, d\mu(x) \leq \left( \int_B |T(f\chi_{2B})(x)|^2 \, d\mu(x) \right)^{1/2} \frac{\mu(B)^{1/2}}{\mu(2B)}$$
$$\leq C \left( \int |f\chi_{2B}|^2 \, d\mu \right)^{1/2} \frac{\mu(B)^{1/2}}{\mu(2B)}$$
$$\leq C \left( \frac{1}{\mu(2B)} \int_{2B} |f(z) - m_{2B}(f)|^2 \, d\mu(z) \right)^{1/2}$$
$$+ C \left( \frac{1}{\mu(2B)} \int_{2B} |m_{2B}(f)|^2 \, d\mu(z) \right)^{1/2} \leq C \|f\|_{\mathcal{L}ip(\alpha,\mu)} r^\alpha.$$

For $A_2$, note first that in case $B = U$, $A_2$ doesn't appear. Otherwise, since $x \notin \text{supp}(f\chi_{2U\setminus 2B})$ and $\mu(2U) \leq C(4r)^n$, we have

$$\frac{1}{\mu(2B)} \int_B \int_{2U\setminus 2B} |k(x,z)| \, |f(z)| \, d\mu(z) \, d\mu(x)$$
$$\leq \frac{1}{\mu(2B)} \int_B \frac{C}{r^n} \int_{2U} |f(z)| \, d\mu(z) \, d\mu(x)$$
$$\leq C \frac{\mu(2U)}{r^n} \frac{1}{\mu(2U)} \int_{2U} |f(z) - m_{2U}(f)| \, d\mu(z) \leq C \|f\|_{\mathcal{L}ip(\alpha,\mu)} r^\alpha.$$

Finally for $A_3$ using the $\epsilon$-regularity of the kernel, we have

$$A_3 \leq \frac{1}{\mu(2B)} \int_B \int_{\mathbb{R}^d\setminus 2U} \frac{|x - x_U|^\epsilon}{|x_U - z|^{n+\epsilon}} |f(z)| \, d\mu(z) \, d\mu(x)$$
$$\leq C \frac{\mu(B)}{\mu(2B)} \|f\|_{\mathcal{L}ip(\alpha,\mu)} r^\alpha \leq C \|f\|_{\mathcal{L}ip(\alpha,\mu)} r^\alpha.$$

The inequalities above are justified by the fact that, for $z \notin 2U$

$$|f(z)| = |f(z) - m_{2U}(f)| \leq \frac{1}{\mu(2U)} \int_{2U} |f(z) - f(y)|\,\mathrm{d}\mu(y)$$
$$\leq (2\,|z - x_U|)^{\alpha}\,\|f\|_{\mathcal{L}ip(\alpha,\mu}$$

and also, that

$$\int_{\mathbb{R}^d \setminus 2U} \frac{1}{|z - x_U|^{n+\epsilon-\alpha}}\,\mathrm{d}\mu(z) \leq Cr^{\alpha-\epsilon}$$

(see [GG]). This concludes the proof of the theorem.

□

DPTO. DE MATEMÁTICAS, C-XV, UNIVERSIDAD AUTÓNOMA, 28049, MADRID, SPAIN
  *E-mail address*: `jose.garcia-cuerva@uam.es`

DEPARTMENT OF MATHEMATICS, DEPAUL UNIVERSITY, CHICAGO, ILLINOIS, 60614, U. S. A.
  *E-mail address*: `aegatto@condor.depaul.edu`